\documentclass[aop,seceqn,MSNbibl,citesort,dvips]{arximspdf}
\usepackage{mathrsfs}

\doi{10.1214/10-AOP548}
\volume{39}
\issue{1}
\pubyear{2011}
\firstpage{224}
\lastpage{251}

\makeatletter

\newproclaim{definition}{Definition}[section]
\newtheorem{lemma}[definition]{Lemma}
\newtheorem{theorem}[definition]{Theorem}
\newproclaim{remarks}[definition]{Remarks}
\newtheorem{corollary}[definition]{Corollary}
\newproclaim{remark}[definition]{Remark}
\newtheorem{proposition}[definition]{Proposition}

\newcommand{\eqref}[1]{(\ref{#1})}
\newcommand{\cmss}{\mathsf}
\newcommand{\Lvec}{L_{\operatorname{vec}}}
\newcommand{\supp}{\operatorname{supp}}
\newcommand{\textd}{{d}}
\newcommand{\texti}{{i}}
\newcommand{\texte}{{e}}
\newcommand{\BB}{\mathcal B}
\newcommand{\EE}{\mathcal E}
\newcommand{\GG}{\mathcal G}
\newcommand{\HH}{\mathcal H}
\newcommand{\LL}{\mathcal L}
\newcommand{\NN}{\mathcal N}
\newcommand{\RR}{\mathcal R}
\newcommand{\B}{\mathbb B}
\newcommand{\R}{\mathbb R}
\newcommand{\Z}{\mathbb Z}
\newcommand{\scrE}{\mathscr{E}}
\newcommand{\scrF}{\mathscr{F}}
\newcommand{\twoeqref}[2]{(\ref{#1})--(\ref{#2})}
\newcommand{\cc}{{\mathrm{c}}}
\newcommand{\pt}{p_{\mathrm{t}}}
\newcommand{\ffrac}[2]{{#1}/{#2}}
\newcommand{\hate}{\hat{{e}}}
\newcommand{\frakq}{\mathfrak q}
\newcommand{\Hspace}{\HH}

\makeatother

\begin{document}
\begin{frontmatter}

\title{Scaling limit for a class of gradient fields with nonconvex potentials}
\runtitle{Scaling limit of gradient fields}

\begin{aug}
\author[A]{\fnms{Marek} \snm{Biskup}\corref{}\thanksref{t1}\ead[label=e1]{biskup@math.ucla.edu}\ead
[label=u1,url]{http://www.math.ucla.edu/\texttildelow biskup/}} and
\author[B]{\fnms{Herbert} \snm{Spohn}\ead[label=e2]{spohn@ma.tum.de}\ead[label=u2,url]{http://www-m5.ma.tum.de/pers/spohn/}}
\runauthor{M. Biskup and H. Spohn}
\affiliation{University of California at Los Angeles and University of
South Bohemia, and Technische Universit\"at M\"unchen}
\address[A]{Department of Mathematics\\
University of California at Los Angeles\\
Los Angeles\\
California 90095-1555\\
USA\\
and\\
School of Economics\\
University of South Bohemia\\
Studentsk\'a 13\\
CZ-37005 \v Cesk\'e Bud\v ejovice\\
Czech Republic\\
\printead{e1}\\
\printead{u1}}
\address[B]{Zentrum Mathematik\\
Technische Universit\"at M\"unchen\\
Boltzmannstrasse 3\\
D-85748 Garching bei M\"unchen\\
Germany\\
\printead{e2}\\
\printead{u2}}
\end{aug}

\thankstext{t1}{Supported by NSF Grant DMS-05-05356.}

\received{\smonth{4} \syear{2007}}
\revised{\smonth{8} \syear{2009}}

%
\begin{abstract}
We consider gradient fields $(\phi_x\dvtx x\in\Z^d)$ whose law takes
the Gibbs--Boltz\-mann form $Z^{-1}\exp\{-\sum_{\langle
x,y\rangle}V(\phi_y-\phi_x)\}$, where the sum runs over nearest
neighbors. We assume that the potential $V$ admits the representation
\[
V(\eta):= - \log\int\varrho(\textd\kappa) \exp\biggl[-\frac{1}{2}\kappa\eta^2\biggr],
\]
where $\varrho$ is a positive measure with compact support
in $(0,\infty)$. Hence, the potential $V$ is symmetric, but nonconvex
in general. While for strictly convex $V$'s, the translation-invariant,
ergodic gradient Gibbs measures are completely characterized by their
tilt, a nonconvex potential as above may lead to several ergodic
gradient Gibbs measures with zero tilt. Still, every ergodic, zero-tilt
gradient Gibbs measure for the potential $V$ above scales to a
Gaussian free field.
\end{abstract}

%

\end{frontmatter}

\section{Introduction}

Gradient fields belong to a class of models that arise in equilibrium
statistical mechanics, for example, as approximations of critical
systems and as effective interface models. Although their definition is
rather simple and, in fact, quite a lot is known (see the reviews by
Funaki \cite{Funaki}, Velenik \cite{Velenik} or Sheffield
\cite{Sheffield}), there is still much to be learned. In this note, we study
gradient fields on a lattice. Here, the field is a collection of
real-valued random variables $\phi:=\{\phi_x\dvtx x\in\Z^d\}$ and the
distribution of $\phi$ on $\R^{\Z^d}$ is given by the formal expression
%
\begin{equation}
\label{grad-field}
\frac1Z\exp\biggl\{-\sum_{\langle x,y\rangle}V(\phi_y-\phi_x) \biggr\}\prod_{x\in\Z
^d}\textd\phi_x,
\end{equation}
where $\textd\phi_x$ is the Lebesgue measure, $\langle x,y\rangle$
refers to an unordered nearest-neighbor pair on $\Z^d$ and $V$ is an
even, measurable function, called the \textit{potential}, which is
bounded from below and grows superlinearly at $\pm\infty$.

Of course, to define the measure \eqref{grad-field} precisely, we have
to restrict the above expression to a finite subset of $\Z^d$ and fix
the $\phi$'s on its boundary; $Z$ is then the normalizing constant.
Another way to regularize the expression \eqref{grad-field} is to
consider directly measures on all of $\R^{\Z^d}$ whose conditional
probabilities in finite sets take the above form. In our context, this
standard definition is hampered by the fact that, due to the unbounded
nature of the fields $\phi_x$, no such infinite-volume measures may
exist at all. However, if one restricts attention to (the $\sigma
$-algebra generated by) the gradient variables
%
\begin{equation}
\label{eta-phi}
\eta_{xy}:=\phi_y-\phi_x,\qquad |x-y|=1,
\end{equation}
then infinite-volume measures exist under the above conditions on $V$.
Since the measure depends only on gradients, we refer to such measures
as \textit{gradient Gibbs measures} (GGM), in accordance with Funaki
\cite{Funaki} and Sheffield \cite{Sheffield}.

Throughout, we will focus on translation-invariant GGMs. An important
characteristic is their \textit{tilt}. For a translation-invariant GGM
$\mu$, there exists a unique tilt vector $t\in\R^d$ such that
%
\begin{equation}
E_\mu(\eta_b)=t\cdot b
\end{equation}
for every edge $b$ of $\Z^d$, which we regard as a vector in this
formula. Of course, this definition is only really meaningful for the
GGMs that are ergodic, that is, trivial on the $\sigma$-algebra of
events invariant under all lattice translations. Indeed, in the ergodic
case, $t$ represents the average incline of typical configurations.

For the case of quadratic $V$, the \textit{massless free field}, the
measure \eqref{grad-field} is Gaussian and so many desired
characteristics are amenable to explicit computations. The challenge
for mathematicians has been to develop an equivalent level of
understanding for nonquadratic $V$'s. A good amount of progress in this
direction has been made in the last ten to fifteen years: Brydges and
Yau \cite{BrydgesYau} (and also earlier works, e.g., by Gaw\c edzki
and Kupiainen \cite{GawedzkiKupiainen} and Magnen and S\'en\'eor
\cite{MagnenSeneor}) studied the effect of analytic perturbations of the
quadratic potentials and concluded that the large-scale behavior is
that of the massless free field. Naddaf and Spencer
\cite{NaddafSpencer} proved the same nonperturbatively for strictly convex
potentials $V$ and GGMs with zero tilt. The corresponding extension to
nonzero tilt was obtained by Giacomin, Olla and Spohn
\cite{GiacominOllaSpohn}. For the same class of potentials, Funaki
and Spohn \cite{FunakiSpohn} proved a bijection between the ergodic
GGMs and their tilt. Sheffield \cite{Sheffield} characterized
translation-invariant GGMs by means of a Gibbs variational principle
and extended Funaki and Spohn's results to fields taking only a
discrete set of values. We refer to the reviews by Funaki
\cite{Funaki}, Velenik \cite{Velenik} and Sheffield \cite{Sheffield} for
further results and references.

As a unifying feature, all the (nonperturbative) results mentioned are
based on the strict convexity of the potential $V$, be it for the use
of the Brascamp--Lieb inequality
\cite{NaddafSpencer,GiacominOllaSpohn,FunakiSpohn}, Helffer--Sj\"ostrand
random walk representation \cite{NaddafSpencer,GiacominOllaSpohn},
coupling to Lange\-vin dynamics \cite{FunakiSpohn} and the
cluster-swapping algorithm \cite{Sheffield}.
One would naturally like to have a nonper\-tur\-bative approach that
works even for nonconvex potentials.
With this motivation, Biskup and Koteck\'y \cite{BiskupKotecky}
recently studied the GGMs for nonconvex $V$ that are a log-mixture of
centered Gaussians,
%
\begin{equation}
\label{V-repr}
V(\eta) := -\log\int\varrho(\textd\kappa) \texte^{-1/2\kappa\eta^2},
\end{equation}
where $\varrho$ is a positive measure with compact support in $(0,\infty
)$. Surprisingly, already for the simplest nontrivial case,
%
\begin{equation}
\label{two-deltas}
\varrho:=p\delta_{\kappa_1}+(1-p)\delta_{\kappa_2}
\end{equation}
with $\kappa_1\gg\kappa_2>0$, it was shown that, in $d=2$, there is a
value $p\in(0,1)$ at which one can construct \textit{two} distinct,
translation-invariant, gradient Gibbs measures of zero tilt.

The relevant conclusion from \cite{BiskupKotecky} for the general
theory is that the one-to-one correspondence between ergodic GGMs and
their tilt breaks down once $V$ is sufficiently nonconvex. The next
question which naturally arises is how to understand what happens to
the scaling limit. The purpose of this note is to show that, regardless
of the occurrence of phase transitions, for potentials of the form
\eqref{V-repr}, every translation-invariant, ergodic GGM with zero tilt
scales to a Gaussian free field (GFF).

The proof is based on the fact---utilized already in
\cite{BiskupKotecky}---that \eqref{V-repr} allows us to represent every GGM
as a mixture over Gaussian gradient measures with a random coupling
constant $\kappa_{xy}$ for each edge $\langle x,y\rangle$. Its
covariance is simply the inverse of the operator
%
\begin{equation}
\label{generator}
(\LL_\kappa f)(x):=\sum_{y\dvtx|y-x|=1}\kappa_{xy}[f(y)-f(x)],
\end{equation}
where we take, once and for all, $\kappa_{xy}=\kappa_{yx}$.
The \textit{fluctuations} in the Gaussian measure can be analyzed by
invoking a random walk representation; $\LL_\kappa$ is the generator of
a random walk with symmetric random jump rates, known, equivalently, as
a \textit{random conductance model}. The name arises naturally from the
electrostatic interpretation of this problem (cf. Doyle and Snell
\cite{DoyleSnell}, in which one views $\Z^d$ as a resistor network with
\textit{conductance} $\kappa_{xy}$---or resistivity $1/\kappa
_{xy}$---assigned to an edge $\langle x,y\rangle$). As it turns out
(see Lemma \ref{lemma-ergodic}), if the initial GGM is ergodic, then so
is the law of the conductances. This makes homogenization a possible tool.

Much work has been done in the past two decades on the problem of
random walks with random conductances. For our purposes, it suffices to
invoke two known results: Kipnis and Varadhan's \cite{KipnisVaradhan}
invariance principle (i.e., scaling of the random walk to Brownian
motion) and Delmotte and Deuschel's \cite{DelmotteDeuschel} annealed
derivative heat kernel bounds. (Note that in the Helffer--Sj\"ostrand
random walk representation, as used in
\cite{NaddafSpencer,GiacominOllaSpohn},
one also has to study a random walk in a random environment. However,
this random environment fluctuates in time, while, in our case, it is
static.) This takes care of the fluctuations of the field; an important
technical issue is thus the control of the mean. This is where the
zero-tilt restriction comes in (see Lemma \ref{lemma-Gaussian},
Corollary \ref{cor-tilt} and discussion in Section \ref{sec7}).

\subsection*{Note} While this paper has been in the process of revision,
further developments have occurred in the study of gradient models with
nonconvex potentials. Cotar, Deuschel and M\"uller \cite{CDM} have
shown that for nonconvex perturbations of potentials $V$ where the size
of the nonconvex region is small compared to typical fluctuations of
the field, the conclusions are as in the convex case. (Their precise
condition is a bound on the $L^1$-norm of the negative part of the
second derivative.) This is a high-temperature result; work in progress
by Adams, Koteck\'y and M\"uller \cite{AKM} addresses the
low-temperature case when nonconvexities are allowed only sufficiently
far away from the absolute minimum of $V$. Our contribution remains
valuable despite these advances as it applies to all potentials of the
type \eqref{V-repr}, including those for which phase coexistence occurs.

This paper is organized as follows. In Section \ref{sec2}, we precisely
define the concept of the gradient Gibbs measure and state our main
theorem. In Section \ref{sec3}, we introduce the extended gradient
Gibbs measures and characterize their conditional marginals. This will
naturally lead to the aforementioned connections with random walks in
reversible random environments. To keep the main line of the argument
intact, we first finish proving our main result in Section \ref{sec4}
and only then expound on the random walk connections in Section
\ref{sec5}. Section \ref{sec7} is devoted to the discussion of the
limitations to zero tilt and some open questions concerning gradient
Gibbs measures.

\section{Main results}
\label{sec2}

\subsection{Gradient Gibbs measures}
As mentioned above, infinite-volume measures on the field variables
$(\phi_x)$ may not always exist, particularly in sufficiently low
dimensions. To make our statements uniform in dimension, we will focus
attention on the gradient variables. However, not even that will be
entirely straightforward because the gradient variables satisfy a host
of ``hard-core'' constraints which, in a sense, encapsulate most of the
interaction. Since $\eta$ is gradient, one has
%
\begin{equation}
\label{loop}
\eta_{x_1,x_2}+\eta_{x_2,x_3}+\eta_{x_3,x_4}+\eta_{x_4,x_1}=0,
\end{equation}
whenever $(x_1,\ldots,x_4)$ are the vertices of a cycle in $\Z^d$ of
length four. We will often write $\eta_b$ for the \textit{positively
oriented} edge $b$ in $\Z^d$. Throughout, we will only work with
positively oriented edges and will use $\B(\Lambda)$ to denote the set
of such edges with both endpoints in the set $\Lambda\subset\Z^d$.

The constraints \eqref{loop} are implemented at the level of the a
priori measure which is defined as follows.
Fix a configuration $\eta\in\R^{\B(\Z^d)}$ that obeys \eqref{loop} and,
for $\Lambda\subset\Z^d$ finite, let $\nu_\Lambda^{(\eta_{\B(\Lambda
)^\cc})}$ be the Lebesgue measure on the linear subspace of
configurations $(\eta_b')$ such that $\eta_b'=\eta_b$ for all $b\notin
\B(\Lambda)$ and that $\eta'$ obeys the constraints \eqref{loop}. Note
that if $\bar\phi$ is a configuration such that $\eta_{xy}=\bar\phi
_y-\bar\phi_x$ for every nearest-neighbor pair $\langle x,y\rangle$,
then $\nu_\Lambda^{(\eta_{\B(\Lambda)^\cc})}$ is, to within a
normalization constant, the projection to gradient variables of the
Lebesgue measure on $\{\phi_x\dvtx x\in\Lambda\}$ subject to the
boundary condition $\bar\phi$.

Next, we will give a precise definition of the notion of gradient Gibbs
measure. For a finite $\Lambda\subset\Z^d$, consider the \textit
{specification} $\gamma_\Lambda$, which is a measure in the first
coordinate and a function of the boundary condition in the second
coordinate, defined by
%
\begin{eqnarray}
\label{spec}
&&\gamma_\Lambda\bigl(\textd\eta_{\B(\Lambda)}|\eta_{\B(\Lambda)^\cc}\bigr)
\nonumber\\[-8pt]\\[-8pt]
&&\qquad:=\frac1{Z_\Lambda(\eta_{\B(\Lambda)^\cc})}
\exp\biggl\{ - \mathop{\sum_{\langle x,y\rangle}}_{x\in\Lambda, y\in\Lambda
\cup\partial\Lambda}
V(\eta_{xy}) \biggr\}
\nu_\Lambda^{(\eta_{\B(\Lambda)^\cc})}\bigl(\textd\eta_{\B(\Lambda
)}\bigr).\nonumber
\end{eqnarray}
Here, $Z_\Lambda(\eta_{\B(\Lambda)^\cc})$ is the normalizing constant.
\begin{definition}
Let $\scrE_{\B(\Lambda)}:=\sigma(\{\eta_b\dvtx b\in\B(\Lambda)\})$. We
say that a measure $\mu$ on $\R^{\B(\Z^d)}$ is a \textit{gradient Gibbs
measure} if the regular conditional probability $\mu(-|\scrE_{\B(\Lambda
)^\cc})$ in any finite $\Lambda\subset\Z^d$ satisfies
%
\begin{equation}
\mu\bigl(- |\scrE_{\B(\Lambda)^\cc}\bigr)(\eta)=\gamma_\Lambda\bigl(-|\eta_{\B(\Lambda
)^\cc}\bigr)
\end{equation}
for $\mu$-a.e. $\eta$.
\end{definition}

Most of this paper is restricted to translation-invariant gradient
Gibbs measures. To define the required notation, for each $x\in\Z^d$,
let $\tau_x\dvtx\R^{\B(\Z^d)}\to\R^{\B(\Z^d)}$ be the ``translation by
$x$'' which acts on configurations $\eta$ by shifting the origin to
position~$x$,
%
\begin{equation}
(\tau_x\eta)_{yz}:=\eta_{y+x,z+x},\qquad (y,z)\in\B(\Z^d).
\end{equation}
We say that $\mu$ is \textit{translation-invariant} if $\mu\circ\tau
_x^{-1}=\mu$ for all $x\in\Z^d$ and that it is \textit{ergodic} if $\mu
(A)\in\{0,1\}$ for every event $A$ such that $\tau_x^{-1}(A)=A$ for all
$x\in\Z^d$.

\subsection{Scaling limit}
As is usual for problems involving random fields, we will interpret
samples from gradient Gibbs measures as random linear functionals on an
appropriate space of functions. Let $C^\infty_0(\R^d)$ denote the set
of all infinitely differentiable functions $f\dvtx\R^d\to\R$ with
compact support.
Given a configuration $\eta=(\eta_b)$ of gradients satisfying the
conditions \eqref{loop}, we can find a configuration of the field $\phi
=(\phi_x)$ such that \eqref{eta-phi} holds for every nearest-neighbor
pair of sites. The configuration $\phi$ is determined uniquely once we
fix the value at one site, for example,~$\phi_0$. For any function
$f\in C^\infty_0(\R^d)$, we introduce the random linear functional
%
\begin{equation}
\phi(f):=\int\textd x\, f(x) \phi_{\lfloor x\rfloor},
\end{equation}
which, under the condition
%
\begin{equation}
\label{zero-sum}
\int\textd x\, f(x)=0,
\end{equation}
does not depend on the choice of the special value $\phi_0$.

The functional $\phi(f)$ can be naturally extended to a somewhat larger
space, defined as follows. Let $\Delta$ denote the Laplace differential
operator in $\R^d$ and consider the set
%
\begin{equation}
\Hspace_0:=\{\Delta g\dvtx g\in C^\infty_0(\R^d)\}.
\end{equation}
Note that each $f\in\Hspace_0$ automatically obeys \eqref{zero-sum}.
The set $\Hspace_0$ is endowed with a natural quadratic form $f\mapsto
(f,f)+(f,-\Delta^{-1}f)$, defined as
%
\begin{equation}
(\Delta g,\Delta g)+(\Delta g,-\Delta^{-1} \Delta g)=\int_{\R^d}\textd
x \bigl(|\Delta g(x)|^2+|\nabla g(x)|^2\bigr).
\end{equation}
We thus define the norm
%
\begin{equation}
\Vert f\Vert_{\Hspace} := [(f,f)+(f,-\Delta^{-1} f)]^{1/2}
\end{equation}
and let $\Hspace$ be the completion of $\Hspace_0$ in this norm.
Note that $\Hspace$ corresponds to the case $k=-\ffrac12$ in the family
of Sobolev spaces $W^{k,2}(\R^d)$. The condition that $(f,-\Delta^{-1}
f)<\infty$ is natural once we realize that this quantity will represent
the variance of the limiting Gaussian field.

The extension of $\phi$ to $\Hspace$ is implied by the following lemma.
\begin{lemma}
\label{lemma-cont}
Suppose that $\varrho$ in \eqref{V-repr} has support bounded away from
zero and let $\mu$ be a translation-invariant, ergodic, zero-tilt
gradient Gibbs measure for the potential $V$. There then exists a
constant $c<\infty$ such that for each $f\in\Hspace_0$,
%
\begin{equation}
\label{4.1}
\Vert\phi(f)\Vert_{L^2(\mu)}\le c\Vert f\Vert_{\Hspace}.
\end{equation}
In particular, $\phi$ extends to a linear functional $\phi\dvtx\Hspace
\to\R$.
\end{lemma}

Note that \eqref{4.1} means that the map $f\mapsto\phi(f)$ is
continuous in $L^2$-norm. If we want to avoid questions about
accumulations of null sets, this permits us to work with only a
countable number of $f$'s at any each time. [In particular, we do not
claim that $f\mapsto\phi(f)$ is continuous in any pointwise sense.]
This will not pose any problems because we will content ourselves with
the following (weaker) definition of a Gaussian free field based on the
standard approach via Gaussian Hilbert spaces (cf. Sheffield
\cite{SheffieldGFF}, Section 2.4).
\begin{definition}
\label{def-GFF}
We say that a family $\{\psi(f)\dvtx f\in\Hspace\}$ of random variables
on a probability space $(\Omega,\scrF,P)$ is a Gaussian free field if
the map $f\mapsto\psi(f)$ is linear a.s. and each $\psi(f)$ is Gaussian
with mean zero and variance
%
\begin{equation}
E(\psi(f)^2)=(f,-\Delta^{-1}f).
\end{equation}
\end{definition}

Our goal is to show that the family of random variables $\{\phi(f)\dvtx
f\in\Hspace\}$ has, asymptotically, in the scaling limit, the law of a
linear transformation of a Gaussian free field. To pass to this limit,
we have to impose the condition that the test functions are slowly
varying, which we take to be on the scale $\varepsilon^{-1}$.
For $\varepsilon>0$ and a function $f\dvtx\R^d\to\R$, let
%
\begin{equation}
f_\varepsilon(x):=\varepsilon^{(d/2+1)}f(\varepsilon x)
\end{equation}
and note that the normalization ensures that
%
\begin{eqnarray}
\Vert f_\varepsilon\Vert_{\Hspace}^2 &=& (f_\varepsilon,f_\varepsilon)+( f_\varepsilon
,(-\Delta)^{-1}f_\varepsilon)
\nonumber\\[-8pt]\\[-8pt]
&=&\varepsilon^2(f,f)+(f,(-\Delta)^{-1}f)\mathop{\le}_{\varepsilon\le1}\Vert
f\Vert_{\Hspace}^2.\nonumber
\end{eqnarray}
Let $\phi_\varepsilon$ denote the linear functional acting on $f\in
C^\infty_0(\R^d)$ via
%
\begin{equation}
\label{phi_epsilon}
\phi_\varepsilon(f) :=\phi(f_\varepsilon)=\int\textd x\, f(x)
\bigl(\varepsilon^{(-d/2 +1)}
\phi_{\lfloor x/\varepsilon\rfloor}\bigr).
\end{equation}
The main theorem is the Gaussian scaling limit for $\phi_\varepsilon(f)$.
\begin{theorem}[(Scaling to GFF)]
\label{thm-main}
Suppose that $V$ is as in \eqref{V-repr} with $\varrho$ compactly
supported in $(0,\infty)$. Let $\mu$ be a gradient Gibbs measure for
the potential $V$ which we assume to be ergodic with respect to the
translations of $\Z^d$ and to have zero tilt. Then, for every $f\in
\Hspace$,
%
\begin{equation}
\label{limit}
\lim_{\varepsilon\downarrow0} E_\mu\bigl(\texte^{\texti\phi_\varepsilon(f)}\bigr) =
\exp\biggl\{\frac12\int\textd x\, f(x) (Q^{-1}f)(x) \biggr\},
\end{equation}
where $Q^{-1}$ is the inverse of the operator
%
\begin{equation}
\label{Q-operator}
Qf:=\sum_{i,j=1}^d q_{ij}\frac{\partial^2}{\partial x_j\,\partial x_i}f,
\end{equation}
with $(q_{ij})$ denoting some positive semidefinite, nondegenerate,
$d\times d$ matrix. In other words, the law of $\phi_\varepsilon$ on the
linear dual $\EE'$ of any finite-dimensional linear subspace $\EE\subset
\Hspace$ converges weakly to that of a Gaussian field with mean zero
and covariance $(-Q)^{-1}$.
\end{theorem}
\begin{remarks}
There follow some additional observations and remarks concerning the
model under consideration and the results above.
\begin{enumerate}[(2)]
\item[(1)] Since $(-Q)$ is dominated by a multiple of $(-\Delta)$ from
below, the integral in \eqref{limit}, inter\-pre\-ted as the quadratic
form $(f,Q^{-1}f)$, is well defined for all $f\in\Hspace$.
\item[(2)]
Note that in \eqref{phi_epsilon}, the individual $\phi$'s get scaled by
$\varepsilon^{(-\ffrac d2+1)}$, not $\varepsilon^{-\ffrac d2}$ as one might
expect from the conventional central limiting reasoning. This has to do
with the fact that the variables $(\phi_x)$ are strongly correlated.
These correlations are weaker for the gradients
$\eta_{xy}:=\phi_y-\phi_x$ which adhere to the ``usual'' central limit
scaling. In $d=1$ and for general potentials $V$, the increments
$\eta_b$ are in fact i.i.d. and the scaling limit follows from the
standard central limit theorem.
\item[(3)] In $d>1$, the matrix
$(q_{ij})$ is not necessarily a multiple of unity since, in general,
$\mu$ is not guaranteed to be invariant under reflections and rotations
of~$\Z^d$. [Nevertheless, we expect that every zero-tilt,
translation-invariant, ergodic measure for the isotropic interaction
\eqref{spec} will inherit these symmetries.] To get convergence of
$\phi_\varepsilon$ to GFF in the sense of Definition \ref{def-GFF}, one
must thus scale the argument of $\phi$ by the root of the
corresponding eigenvalue of $\mathfrak{q}$ in each of its principal
directions.
\item[(4)] The absence of strict convexity does not permit
us to use the general argument of Funaki and Spohn \cite{FunakiSpohn}
for the \textit{existence} of an ergodic GGM with zero (or any other
prescribed) tilt. To show that such GGMs do exist---and that our
Theorem \ref{thm-main} is not vacuous---we note that, by Lemma 4.8 of
Biskup and Koteck\'y~\cite{BiskupKotecky}, every weak limit of torus
measures exhibits exponential concentration of the empirical tilt; one
then just needs to choose any ergodic component. Note that this lemma
applies only to zero tilt (cf. \cite{BiskupKotecky}, Remark 4.9).
\item[(5)] The restriction to zero tilt is actually a significant
drawback of our analysis. The main reason is our inability to
characterize the scaling limit of the so-called corrector for the
corresponding random walk problem. See Section \ref{sec7} for more
details.
\item[(6)] In the example studied by Biskup and Koteck\'y
\cite{BiskupKotecky} [cf. \eqref{two-deltas}], the two GGMs coexisting
at the transitional value $\pt$ of $p$ were proven to exhibit different
characteristic fluctuations. It follows that the corresponding scaling
limits will be distinguished by their stiffness coefficients $q_{ij}$.

Moreover, by Theorem 2.5 of \cite{BiskupKotecky}, for
$\kappa_1\gg\kappa_2$, the transition in the $d=2$ model with
\eqref{two-deltas} lies on a self-dual line, that is,
%
\begin{equation}
\frac{\pt}{1-\pt}= \biggl(\frac{\kappa_2}{\kappa_1} \biggr)^{\ffrac14}.
\end{equation}
The transition presumably stays on this line even as one slides the
ratio $\kappa_1/\kappa_2$ toward one. However, it disappears before
$\kappa_1/\kappa_2$ hits one because, for $\kappa_1\approx\kappa_2$,
the potential $V$ is convex and so there is only one GGM with zero tilt
\cite{FunakiSpohn}. At such a point of disappearance, physicists often
expect nontrivial critical fluctuations. Notwithstanding, our results
show that this is not the case.

\item[(7)]
We avoid the context of the ``stronger'' definition of GFF as a random
element in an appropriate Banach space (cf. Gross \cite{Gross} or
Sheffield \cite{SheffieldGFF}, Section 2.2). This definition is
appealing in $d=1$, where the limiting functional $f\mapsto\psi(f)$
actually admits the integral representation
%
\begin{equation}
\psi(f)=\int_\R f(t)\psi_t\,\textd t
\end{equation}
with $t\mapsto\psi_t$ denoting a continuous diffusion with generator
$Q$, but in $d>1$, the corresponding field becomes less and less
regular with increasing dimension and the appeal is lost. However, this
context would be ideal if one wished to discuss the notion of tightness
and convergence in law for the limit in Theorem \ref{thm-main}.
\end{enumerate}
\end{remarks}

Both Lemma \ref{lemma-cont} and Theorem \ref{thm-main} are proved in
Section \ref{sec4}.

\section{Extended gradient Gibbs measures}
\label{sec3}

\subsection{Coupling to random conductance model}
The key idea underlying the representation \eqref{V-repr} is that the
auxiliary variable $\kappa$ in the expression for $V$ may be elevated
to a genuine degree of freedom associated with the corresponding edge.
Specifically, given a gradient Gibbs measure $\mu$ with potential \eqref
{V-repr}, for each finite $\Lambda\subset\B(\Z^d)$, consider the
measure $\tilde\mu_\Lambda$ on $\R^{\B(\Z^d)}\times\R^\Lambda$ defined by
%
\begin{equation}
\label{consistent}
\tilde\mu_\Lambda({\mathcal A}\times\BB):=\int_\BB\prod_{b\in\Lambda}\varrho
(\textd\kappa_b) E_\mu\biggl( \mathsf{1}_{\mathcal A}\prod_{b\in\Lambda}\texte^{V(\eta
_b)-1/2\kappa_b\eta_b^2} \biggr),
\end{equation}
where ${\mathcal A}\subset\R^{\B(\Z^d)}$ and $\BB\subset\R^\Lambda$ are Borel
sets. The representation \eqref{V-repr} ensures that $(\tilde\mu_\Lambda
)$ is a consistent family of measures; by Kolmogorov's extension
theorem, these are projections from a unique measure $\tilde\mu$ onto
configurations $(\eta_b,\kappa_b)\in\R^{\B(\Z^d)}\times\R^{\B(\Z^d)}$.
The restriction of $\tilde\mu$ to the $\eta$'s gives us back $\mu$; we
call $\tilde\mu$ an \textit{extension} of $\mu$. The measure $\tilde\mu$
is Gibbs for the Hamiltonian $\sum_{\langle x,y\rangle}\frac12\kappa
_{xy}\eta_{xy}^2$, so we will refer to it as an \textit{extended gradient
Gibbs measure} (see Biskup and Koteck\'y \cite{BiskupKotecky} for
further facts on extended GGMs).

To ease the notation, whenever $b$ is an edge between $x$ and $y$, we
may interchangeably write $\kappa_b$ and $\kappa_{xy}$ for the same
quantity. Furthermore, for the same reasons, it will even be convenient
to assume
that
%
\begin{equation}
\kappa_{xy}=\kappa_{yx},\qquad |x-y|=1
\end{equation}
and work with the $\kappa$'s as symmetric objects.

We proceed with a series of lemmas that characterize the properties of
$\tilde\mu$.
\begin{lemma}
\label{lemma-iid}
Let $\mu$ be a gradient Gibbs measure for the potential $V$ and let
$\tilde\mu$ be its extension to $\R^{\B(\Z^d)}\times\R^{\B(\Z^d)}$.
Consider the $\sigma$-field $\scrE:=\sigma(\{\eta_b\dvtx b\in\B(\Z^d)\}
)$. For $\tilde\mu$-a.e. $\eta$, the regular conditional distribution
$\tilde\mu(-|\scrE)(\eta)$, regarded as a measure on the $\kappa$'s,
takes the product form
%
\begin{equation}
\tilde\mu(\textd\kappa|\scrE)(\eta)=\bigotimes_{b\in\B(\Z^d)}\bigl[\texte
^{V(\eta_b)-1/2\kappa_b\eta_b^2}\varrho(\textd\kappa_b)\bigr].
\end{equation}
\end{lemma}
\begin{pf}
Recall that $\scrE_{\B(\Lambda)}:=\sigma(\{\eta_b\dvtx b\in\B(\Lambda)\}
)$. The identity \eqref{consistent} implies that $\tilde\mu_\Lambda$
coincides with $\tilde\mu$ on $\R^{\B(\Lambda)}\times\R^{\B(\Lambda)}$.
However, $\mu_\Lambda(-|\scrE_\Lambda)$ has the desired product form by
definition and so the claim follows by standard approximation arguments.
\end{pf}
\begin{lemma}
\label{lemma-ergodic}
Let $\mu$ be a gradient Gibbs measure and let $\tilde\mu$ be its
extension to $\R^{\B(\Z^d)}\times\R^{\B(\Z^d)}$.
If $\mu$ is translation-invariant and ergodic, then so is $\tilde\mu$.
\end{lemma}
\begin{pf}
The uniqueness of the extension of measures \eqref{consistent} implies
that $\tilde\mu$ is translation-invariant if $\mu$ is
translation-invariant and so it remains to prove that ergodicity is
also inherited. Let $A\subset\R^{\B(\Z^d)}\times\R^{\B(\Z^d)}$ be a
translation-invariant event, that is, $(\eta,\kappa)\in A$ if and only
if $(\tau_x\eta,\tau_x\kappa)\in A$ for all $x$. Our task is to show
that $\tilde\mu(A)\in\{0,1\}$.

First, we invoke the ergodicity of $\mu$. Consider the function
%
\begin{equation}
f(\eta):=E_{\tilde\mu}(\mathsf{1}_A|\scrE)(\eta).
\end{equation}
Since $A$ and $\tilde\mu$ are translation-invariant, we have
%
\begin{equation}\quad
f(\tau_x\eta)=E_{\tilde\mu}(\mathsf{1}_A|\scrE)(\tau_x\eta)=E_{\tilde\mu}(\mathsf{1}
_A\circ\tau_x^{-1}|\scrE)(\eta)=f(\eta),\qquad
\tilde\mu\mbox{-a.s.}
\end{equation}
But $f$ is $\scrE$-measurable and the restriction of $\tilde\mu$ to
$\scrE$ is $\mu$, which we assumed to be ergodic. Hence, $f$ is
constant almost surely. Let $c$ denote this constant.

We will use a standard approximation argument to show that $c\in\{0,1\}
$. Since $A$ is an event from the product $\sigma$-algebra, there
exists a sequence of events
%
\begin{equation}
A_n\in\sigma(\{\eta_{xy},\kappa_{xy}\dvtx|x-y|=1, |x|\le n\})
\end{equation}
such that
%
\begin{equation}
\tilde\mu(A\triangle A_n) \mathop{\longrightarrow}_{n\to\infty}0.
\end{equation}
The bound
%
\begin{equation}
\Vert\mathsf{1}_{A_n}-\mathsf{1}_A\Vert_{L^1(\tilde\mu)}\le\tilde\mu(A\triangle A_n)
\end{equation}
then shows that $\mathsf{1}_{A_n}\to\mathsf{1}_A$ in $L^1(\tilde\mu)$. Since $A$ is
translation-invariant, we have $\mathsf{1}_A=\mathsf{1}_A\mathsf{1}_{\tau_x^{-1}(A)}$. Each
indicator can\vspace*{2pt} be approximated by the indicator of the event $A_n$; a
simple bound gives
%
\begin{equation}
\bigl\Vert\mathsf{1}_{A_n}\mathsf{1}_{\tau_x^{-1}(A_n)}-\mathsf{1}_A\mathsf{1}_{\tau_x^{-1}(A)}\bigr\Vert
_{L^1(\tilde\mu)}\le2\tilde\mu(A\triangle A_n).
\end{equation}
For $x$ with $|x|>2n+1$, the fact that $\tilde\mu(-|\scrE)$ is a
product measure (cf. Lem\-ma~\ref{lemma-iid}) implies that $A_n$ and $\tau
_x^{-1}(A_n)$ are independent. Hence,
%
\begin{equation}
E_{\tilde\mu}\bigl(\mathsf{1}_{A_n}\mathsf{1}_{\tau_x^{-1}(A_n)}|\scrE\bigr)=E_{\tilde\mu}(\mathsf{1}
_{A_n}|\scrE)E_{\tilde\mu}\bigl(\mathsf{1}_{\tau_x^{-1}(A_n)}|\scrE\bigr).
\end{equation}
Rolling the approximations backward, we then conclude that the
left-hand side converges to $c$ in $L^1(\tilde\mu)$, while the
right-hand side converges to $c^2$ (note that all expectations are
bounded). It follows that $c=c^2$, that is, $c\in\{0,1\}$. As
%
\begin{equation}
\tilde\mu(A)=E_{\tilde\mu}(f)=c,
\end{equation}
the proof is finished.
\end{pf}

\subsection{Random walk connections}
Our next goal will be to characterize also the conditional measure
given the $\kappa$'s. This will, in turn, require some facts from the
theory of random walks with random conductances. We will frequently
borrow facts from an associated potential theory which will be
expounded in Section \ref{sec5}.

Let us choose a configuration $\kappa=(\kappa_b)$ with $\kappa_b\in
(0,\infty)$ and recall the formula \eqref{generator} for the generator
$\LL_\kappa$ of the random walk among conductances $\kappa$. We will
focus on the action of $\LL_\kappa$ on functions of both the
environment $\kappa$ and the position $x$ that satisfy the following
\textit{shift covariance} property:
%
\begin{equation}
\label{shift-invariance}
g(\kappa,x+b)-g(\kappa,x)=g(\tau_{x}\kappa,b),
\end{equation}
with $x\in\Z^d$ and $b$ a coordinate unit vector in $\R^d$,
subject to the condition
%
\begin{equation}
\label{zero-cond}
g(\kappa,0)=0.
\end{equation}
This makes the function completely determined by its values at the
neighbors of the origin. A function of this kind is said to be \textit
{harmonic} for the above random walk if
%
\begin{equation}
\label{harmonic}
\LL_\kappa g(\kappa,\cdot)=0
\end{equation}
for (almost) every $\kappa$. As it turns out, harmonic, shift-covariant
functions are uniquely determined (a.s.) by their mean with respect to
ergodic measures on the conductances.
\begin{lemma}
\label{lemma-harmonic}
Let $\nu$ be a translation-invariant, ergodic probability measure on
configurations $\kappa=(\kappa_b)\in\R^{\B(\Z^d)}$ such that $\nu
(\varepsilon\le\kappa_b\le\ffrac1\varepsilon)=1$ for some $\varepsilon>0$. Let
$g\dvtx\R^{\B(\Z^d)}\times\Z^d\to\R$ be a measurable function which is:
\begin{enumerate}
\item[(1)]
harmonic in the sense of \eqref{harmonic}, $\nu$-a.s.;
\item[(2)]
shift-covariant in the sense of (\ref{shift-invariance}) and (\ref{zero-cond}),
$\nu$-a.s.;
\item[(3)]
square integrable in the sense that $E_\nu|g(\cdot,x)|^2<\infty$ for
all $x$ with $|x|=1$.
\end{enumerate}
If $E_\nu(g(\cdot,x))=0$ for all $x$ with $|x|=1$, then $g(\cdot,x)=0$
a.s. for all $x\in\Z^d$.
\end{lemma}

We defer the proof, and further discussion of the consequences of shift
covariance and harmonicity, to Section \ref{sec5}.
Returning to the gradient fields, we now characterize the conditional
law given the $\kappa$'s.
\begin{lemma}
\label{lemma-Gaussian}
Let $\mu$ be a translation-invariant, ergodic gradient Gibbs measure
with zero tilt and let $\tilde\mu$ be its extension to $\R^{\B(\Z
^d)}\times\R^{\B(\Z^d)}$. Consider the $\sigma$-field $\scrF:=\sigma(\{
\kappa_b\dvtx b\in\B(\Z^d)\})$. For $\tilde\mu$-a.e. $\kappa$, the
conditional law $\tilde\mu(-|\scrF)(\kappa)$, regarded as a measure on
the set of configurations $\{(\phi_x)\in\R^{\Z^d}\dvtx\phi_0=0\}$ with
the $\phi$'s defined from the $\eta$'s via \eqref{eta-phi}, is Gaussian
with mean zero,
%
\begin{equation}
\label{mean-zero}
E_{\tilde\mu}(\phi_x|\scrF)(\kappa)=0,\qquad x\in\Z^d,
\end{equation}
and covariance given by $(-\LL_\kappa)^{-1}$. Explicitly, for each
$f\dvtx\Z^d\to\R$ with finite support and $\sum_xf(x)=0$,
%
\begin{equation}
\label{variance}
\operatorname{Var}_{\tilde\mu} \biggl(\sum_xf(x)\phi_x \big|\scrF\biggr)(\kappa)=\sum
_{x}f(x)(-\LL_\kappa^{-1}f)(x).
\end{equation}
\end{lemma}
\begin{pf}
The fact that the conditional measure is a multivariate Gaussian law
with covariance $\LL_\kappa^{-1}$ can be checked by direct inspection
of \eqref{consistent}. The only nontrivial task is to identify the
mean. First, we note that the loop conditions \eqref{loop} ensure that
there exists a function $u\dvtx\R^{\B(\Z^d)}\times\Z^d\to\R$ such that
%
\begin{equation}
u(\kappa,0)=0
\end{equation}
and
%
\begin{equation}
u(\kappa,x+b)-u(\kappa,x)=E_{\tilde\mu}(\eta_{x,x+b}|\scrF)(\kappa)
\end{equation}
for all unit vectors $b$ in the coordinate directions. We claim that
$u$ is harmonic in the sense of \eqref{harmonic}. Indeed,
%
\begin{equation}
\LL_\kappa u(\kappa,x)=E_{\tilde\mu} \biggl( \sum_{y\dvtx|y-x|=1}\kappa
_{xy}(\phi_y-\phi_x) \big|\scrF\biggr)(\kappa),
\end{equation}
where we write, thanks to the loop conditions, $\eta_{xy}=\phi_y-\phi_x$.
Using the fact that $\tilde\mu$ is Gibbs, we can now also condition on
$\sigma(\phi_y\dvtx y\ne x)$; the conditional measure $\mu_{\{x\}}$ is
Gaussian with the explicit form
%
\begin{equation}\qquad
\mu_{\{x\}}(\textd\phi_x)=\frac1Z\exp\biggl\{-\frac12\phi_x^2\sum_{y\dvtx
|y-x|=1}\kappa_{xy}+\phi_x \sum_{y\dvtx|x-y|=1}\kappa_{xy}\phi_y \biggr\}
\,\textd\phi_x,
\end{equation}
where $Z$ is an appropriate normalization constant. It is easy to check
that the mean of $\phi_x\sum_{y\dvtx|y-x|=1}\kappa_{xy}$ under $\mu_{\{
x\}}$ is exactly $\sum_{y\dvtx|y-x|=1}\kappa_{xy}\phi_y$, proving that
$\LL_\kappa u(\kappa,x)=0$.

Next, we observe that the translation invariance of $\tilde\mu$ implies
that
%
\begin{eqnarray}
u(\tau_{x}\kappa,b)-u(\tau_{x}\kappa,0)
&=& E_{\tilde\mu}(\eta_{0,b}|\scrF)(\tau_{x}\kappa)
\nonumber\\
&=&E_{\tilde\mu}(\eta_{x,x+b}|\scrF)(\kappa)\\
&=&u(\kappa,x+b)-u(\kappa,x)
\nonumber
\end{eqnarray}
and so $u$ is shift-covariant, as defined in
(\ref{shift-invariance}) and (\ref{zero-cond}). Finally, the definition of $u$ and the
fact that $\tilde\mu$ has zero tilt imply that
%
\begin{equation}
\label{u-exp}
E_{\tilde\mu}(u(\cdot,x))=E_{\tilde\mu}(\eta_{0,x})=0,\qquad |x|=1.
\end{equation}
As $u$ obeys all conditions of Lemma \ref{lemma-harmonic}, we have
$E_{\tilde\mu}(\phi_x|\scrF)=u(\cdot,x)=0$ $\tilde\mu$-a.s.
\end{pf}

Our reference to the random walk with generator $\LL_\kappa$ is not
limited to Lem\-ma~\ref{lemma-harmonic}; we will also need to know some
specific properties of this random walk. First, we will need to know
that the position of the walk satisfies a central limit theorem. Let
$X=(X_t)$ denote the continuous-time random walk with the generator $\LL
_\kappa$ and let $P_\kappa^x$ denote the law of $X$ subject to the
initial condition $P_\kappa^x(X_0=x)=1$. The following lemma goes back
to Kipnis and Varadhan \cite{KipnisVaradhan}.
\begin{lemma}[(Annealed central limit theorem)]
\label{lemma-BM}
Let $\mu$ be a measure on $\R^{\B(\Z^d)}$ which is
tran\-sla\-tion-invariant, ergodic and obeys $\mu(\varepsilon\le\omega_b\le\ffrac
1\varepsilon)=1$ for some $\varepsilon>0$. There then exists a positive
semidefinite, nondegenerate, $d\times d$ matrix $\frakq$ such that for
every $t>0$, the annealed distribution $E_\mu P_\kappa^0(\varepsilon
X_{t\varepsilon^{-2}}\in\cdot)$ converges weakly to the law of the
multivariate normal $\NN(0,t\frakq)$.
\end{lemma}

The main result of \cite{KipnisVaradhan} actually shows that the
annealed law of the entire \textit{path} $t\mapsto\varepsilon X_{t\varepsilon
^{-2}}$ converges to that of (a linear transform of) Brownian motion.
However, the above is all that will be needed for the purposes of the
present paper.

Apart from a central limit asymptotics, we will also need an estimate
on the heat kernel of the above random walk. The following lemma is a
consequence of the main result of Delmotte and Deuschel
\cite{DelmotteDeuschel}.
\begin{lemma}[(Heat kernel upper bound)]
Let $\mu$ be a law of the conductances satisfying the ellipticity
condition $\mu(\varepsilon<\kappa_b<\ffrac1\varepsilon)=1$ for some $\varepsilon
>0$. There is then a $c_1<\infty$ such that
%
\begin{equation}
\label{HCUB}
E_\mu|\nabla_i\nabla_jP_\kappa^0(X_t=\cdot)|\le\frac{c_1}{t^{
d/2+1}},\qquad 1\le i,j\le d, t>0.
\end{equation}
Here, $\nabla_i$ is the discrete spatial derivative in the $i$th
coordinate direction, that is, $\nabla_i f(x):=f(x+\hate_i)-f(x)$.
\end{lemma}
\begin{pf}
By formula (1.5b) in \cite{DelmotteDeuschel}, Theorem 1.1,
%
\begin{equation}
E_\mu|\nabla_i\nabla_jP_\kappa^0(X_t=x)|
\le c_1'\frac{\cmss P^{c_2't}(0,x)}{ t},
\end{equation}
where $\cmss P^{c_2't}(0,x)$ is the probability of the continuous-time
simple random walk to be at $x$ at time $c_2't$. This probability is
bounded from above by a constant times $t^{-d/2}$.
\end{pf}

\section{Proof of main result}
\label{sec4}

\subsection{Regularity estimates}
The goal of this section is to prove Theorem \ref{thm-main} concerning
the scaling limit of $\phi(f_\varepsilon)$. We begin by proving
$L^2$-continuity of the random functional $f\mapsto\phi(f)$ on $\Hspace
$, as stated in Lemma \ref{lemma-cont}. For convenience of notation,
whenever $\RR$ is an operator on $\ell^2(\Z^d)$, we will extend it to
an operator on $L^2(\R^d)$ via the formula
%
\begin{equation}
(f,\RR f) :=\int\textd x\,\textd y\, f(x)f(y) \RR(\lfloor x\rfloor,\lfloor
y\rfloor),
\end{equation}
where $\RR(x,y)$ is the kernel of $\RR$ in the canonical basis in $\ell
^2(\Z^d)$.
\begin{pf*}{Proof of Lemma \ref{lemma-cont}}
Let $\tilde\mu$ be the extended gradient Gibbs measure corresponding to
$\mu$. Recall the notation $\LL_\kappa$ for the generator of the random
walk among conductances $\kappa=(\kappa_b)$ and let $\LL$ denote the
generator of the simple random walk (i.e., the special case of $\LL
_\kappa$ when all $\kappa_b=1$). Choose $f\in\{\Delta g\dvtx g\in
C_0^\infty(\R^d)\}$. Lemma \ref{lemma-Gaussian} and the fact that $\phi
(f)$ is linear in the $\eta$'s imply
that
%
\begin{equation}
\Vert\phi(f)\Vert_{L^2(\tilde\mu)}^2=E_{\tilde\mu}((f,-\LL_\kappa^{-1} f)),
\end{equation}
where $(f,\LL_\kappa^{-1} f)$ is as defined above.
By assumption on the support of $\varrho$, we know that $\kappa_b\ge a$
$\tilde\mu$-a.s., by which we immediately have the operator inequalities
%
\begin{equation}
(-\LL_\kappa)\ge a(-\LL) \quad\mbox{and}\quad (-\LL_\kappa)^{-1}\le a^{-1}(-\LL)^{-1}.
\end{equation}
Therefore, it suffices to bound the quadratic form associated with the
(homogeneous) discrete Laplacian $\LL$ by the quadratic form defining
the $\Hspace$-space:
%
\begin{equation}
\label{discrete-continuous}
(f,(-\LL)^{-1} f)\le c \Vert f\Vert_\Hspace^2
\end{equation}
for some constant $c<\infty$ and all $f$ in a dense subset of $\Hspace$.

To this end, we pick $f\in\Hspace$ in the Schwartz class and let
%
\begin{equation}
\hat f(k):=(2\pi)^{-d/2}\int f(x)\texte^{\texti k\cdot x}\,\textd x\,
\end{equation}
be its ($L^2$-norm-preserving) Fourier transform. A direct calculation
now yields
%
\begin{equation}
\label{Fourier}\qquad
(f,(-\LL)^{-1} f)=\int_{[-\pi,\pi]^d}\textd k \frac1{(-\hat\LL)(k)} \Biggl|
\mathop{\sum_{k'\dvtx\exists\ell\in\Z^d}}_{k-k'=2\pi\ell}
\hat f(k')\prod_{j=1}^d\frac{1-\texte^{-\texti k_j'}}{\texti k_j'} \Biggr|^2,
\end{equation}
where
%
\begin{equation}
(-\hat\LL)(k):=4\sum_{j=1}^d \sin^2(k_j/2)
\end{equation}
is the generalized eigenvalue of the lattice Laplacian. Introducing
$-\hat\Delta(k)=|k|^2$ to denote the corresponding quantity for the
continuum Laplacian, we invoke the Cauchy--Schwarz inequality for the
sum over $k'$ to get
%
\begin{equation}
\label{CSbound}\quad
(f,(-\LL)^{-1} f)\le\int_{[-\pi,\pi]^d}\textd k\, \biggl\{ \mathop{\sum
_{k'\dvtx\exists\ell\in\Z^d}}_{k-k'=2\pi\ell}
|\hat f(k')|^2\bigl(1-\hat\Delta(k')^{-1}\bigr) \biggr\} c(k),
\end{equation}
where
%
\begin{equation}
\label{ck-def}
c(k):=\mathop{\sum_{k'\dvtx\exists\ell\in\Z^d}}_{k-k'=2\pi\ell}
\frac1{(-\hat\LL)(k)} \frac{-\hat\Delta(k')}{1-\hat\Delta(k')}
\prod_{j=1}^d \biggl(\frac{2\sin(k_j'/2)}{k_j'} \biggr)^2
\end{equation}
is well defined on the set $\BB:=\{k\in[-\pi,\pi]^d\dvtx k_j\ne0,
j=1,\ldots,d\}$ of full Lebesgue measure in $[-\pi,\pi]^d$. We claim that
%
\begin{equation}
\label{c-sup}
c:=\sup_{k\in\BB}c(k)<\infty.
\end{equation}
We will show this by proving that the summand in \eqref{ck-def} is
bounded by a constant times the product $\prod_j(|k_j'|+1)^{-2}$.
Indeed, for the $k'=k$ term, we use the fact that the ratio $\hat\Delta
(k)/\hat\LL(k)$ is bounded throughout $\BB$, and the same for the
ratios $k\mapsto\sin(k_j/2)/k_j$. When $k'\ne k$, we set $i$ to be the
first index $j$ such that $k_j'\ne k_j$ and bound the $4\sin(k_i'/2)^2$
term by $-\hat\LL(k)$. We then bound the ratio of $\hat\Delta(k)$ terms
by unity and the $j$th term in the product by a constant times
$(|k_j'|^2+1)^{-2}$. [The $\sin(k_i'/2)$ term is not needed because
$|k_i'|\ge2\pi$.]

The product $\prod_j(|k_j'|+1)^{-2}$ is summable over $k'\in k+(2\pi\Z
)^d$ uniformly in $k\in\BB$ and so \eqref{c-sup} is proved.
Bounding $c(k)$ by its supremum in \eqref{CSbound}, we can merge the
sum and the integral to get $\Vert f\Vert_\Hspace^2$. The desired bound
\eqref{discrete-continuous} then follows.
\end{pf*}
\begin{remark}
The inclusion of $L^2$-norm of $f$ in $\Vert f\Vert_\Hspace$ is crucial
for the bound \eqref{discrete-continuous}. Indeed, on the basis of
\eqref{Fourier}, it is not hard to construct functions for which the
ratio $(f,-\LL^{-1}f)/(f,-\Delta^{-1}f)$ is arbitrarily large. This is
caused by the fact that the spectrum of $-\Delta$ extends all the way
to infinity, while that of $-\LL$ is bounded.
\end{remark}

The continuity established in Lemma \ref{lemma-cont} allows us to work
only with smooth and compactly supported test functions. We will
nevertheless need one more regularity bound before we can delve into
the proof of our main result.
\begin{lemma}
\label{lemma-reg}
Let $\mu$ be a translation-invariant law on the conductances subject to
the ellipticity condition $\mu(\varepsilon<\kappa_b<\ffrac1\varepsilon)=1$
for some $\varepsilon>0$. There then exists $c<\infty$ such that whenever
$f=\Delta g$ for some $g\in C_0^\infty(\R^d)$,
%
\begin{equation}
\label{reg-bd}
E_\mu(f,\texte^{t\LL_\kappa}f)\le c\Vert\nabla g\Vert_\infty^2 \lambda
(\supp g)^2\frac{1}{t^{d/2+1}},
\end{equation}
where $\lambda(A)$ is the set function on Borel subsets of $\R^d$
defined by
%
\begin{equation}
\lambda(A):=\sum_{x\in\Z^d}\sum_{i=1}^d\int_{R_i}\textd z\,\mathsf{1}_{\{x+z\in
A\}},
\end{equation}
with $R_i$ denoting the set of points in $[0,1]^d$ with vanishing $i$th
coordinate.
\end{lemma}
\begin{pf}
Translation invariance of $\mu$ and a simple integration by parts tells us
that
%
\begin{eqnarray}\quad
&&
E_\mu(\Delta g,\texte^{t\LL_\kappa}\Delta g)
\nonumber\\
&&\qquad=\sum_{x,y\in\Z^d}\int_{[0,1]^d}\textd z\int_{[0,1]^d}\textd z'\, \Delta
g(x+z)\Delta g(y+z')
E_\mu P^x_\kappa(X_t=y)
\nonumber\\[-8pt]\\[-8pt]
&&\qquad=
\sum_{x,y\in\Z^d}\sum_{i,j=1}^d\int_{R_i}\textd z\int_{R_j}\textd z'\,
\partial_ig(x+z)\,\partial_j g(y+z')\nonumber\\
&&\qquad\quad\hspace*{94pt}{}\times
\nabla_i\nabla_j E_\mu P^x_\kappa(X_t=y),\nonumber
\end{eqnarray}
where $\partial_i$ stands for the partial derivative with respect to
the $i$th coordinate. Restricting the integrations and sums so that the
arguments $x+z$ and $y+z'$ are in the support of $g$, bounding the
partial derivatives by $\Vert\nabla g\Vert_\infty$ and applying the
estimate~\eqref{HCUB}, we obtain the desired bound.
\end{pf}

The consequence of Lemma \ref{lemma-reg} that will concern us is as follows.
\begin{corollary}
\label{cor-tightness}
For $\mu$ as in Lemma \ref{lemma-reg} and any $f\in\{\Delta g\dvtx g\in
C^\infty_0(\R^d)\}$,
%
\begin{equation}
\lim_{M\to\infty}\sup_{0<\varepsilon<1}
\int_M^\infty\textd t\, E_\mu\varepsilon^{-2}(f_\varepsilon,\texte^{t\varepsilon
^{-2}\LL_\kappa} f_\varepsilon)=0.
\end{equation}
\end{corollary}
\begin{pf}
Choose $f$ of the form $f=\Delta g$ and note that $f_\varepsilon=\Delta
g^{(\varepsilon)}$, where $g^{(\varepsilon)}(x):=\varepsilon^{
d/2+1}g(x\varepsilon)$. First, we observe
that
%
\begin{equation}
\label{4.13a}
\bigl\Vert\nabla g^{(\varepsilon)}\bigr\Vert_\infty=\varepsilon^{d/2}\Vert\nabla
g\Vert_\infty.
\end{equation}
Next, we note that, since the support of $g$ is the closure of a
nonempty bounded open set, a simple covering argument tells us that
%
\begin{equation}
\varepsilon^{d}\lambda\bigl(\supp g^{(\varepsilon)}\bigr)
\mathop{\longrightarrow}_{\varepsilon\downarrow0} d|{\supp g}|,
\end{equation}
where $|{\supp g}|$ is the Lebesgue measure of $\supp g$. As a
consequence, there exists a constant $C(g)<\infty$ such that
%
\begin{equation}
\label{4.15a}
\lambda\bigl(\supp g^{(\varepsilon)}\bigr)\le C(g)\varepsilon^{-d},\qquad 0<\varepsilon<1.
\end{equation}
Plugging \eqref{4.13a} and \eqref{4.15a} into \eqref{reg-bd}, we get,
for $\varepsilon\in(0,1)$,
%
\begin{equation}
E_\mu\varepsilon^{-2}(f_\varepsilon,\texte^{t\varepsilon^{-2}\LL_\kappa}
f_\varepsilon)
\le c\Vert\nabla g\Vert_\infty^2 C(g)^2\frac{1}{t^{d/2+1}}.
\end{equation}
The functions on the left (indexed by $\varepsilon$) are uniformly
integrable in $t$ in all $d\ge1$.
\end{pf}

\subsection{Scaling limit}
Having dispensed with regularity considerations, we can now proceed to
establish the principal fact underlying the proof of Theorem \ref{thm-main}.
\begin{proposition}
\label{lemma-scaling}
Let $\mu$ be a translation-invariant, ergodic measure on $\kappa=(\kappa
_b)\in\R^{\B(\Z^d)}$ such that $\mu(\delta\le\kappa_b\le\ffrac1\delta
)=1$ for some $\delta>0$. There then exists a positive semidefinite,
nondegenerate $d\times d$ matrix $\frakq=(q_{ij})$ such that
%
\begin{equation}
\label{4.8}
\lim_{\varepsilon\downarrow0}(f_\varepsilon,(-\LL_\kappa)^{-1}f_\varepsilon) =
(f,(-Q)^{-1}f)
\end{equation}
in $\mu$-probability for each $f\in C_0^\infty(\R^d)\cap\Hspace$,
where $Q$ is defined from $\frakq$ by \eqref{Q-operator}.
\end{proposition}

The key to the proof is the following lemma.
\begin{lemma}
For every $t>0$ and any $f\in C_0^\infty(\R^d)\cap\Hspace$,
%
\begin{equation}
\label{4.9}
\Theta_\varepsilon(t) :=\varepsilon^{-2}(f_\varepsilon,\texte^{ t\varepsilon^{-2}\LL
_\kappa}f_\varepsilon)-(f,\texte^{ tQ}f) \mathop{\longrightarrow}_{\varepsilon
\downarrow0} 0\qquad
\mbox{in }L^2(\mu).
\end{equation}
\end{lemma}
\begin{pf}
Let $(X_t)_{t\ge0}$ be the continuous-time random walk with the
generator $\LL_\kappa$ and let $P_\kappa^x$ denote the law of the walk
started from $x$. By Lemma \ref{lemma-BM}, the annealed law of
$\varepsilon X_{t\varepsilon^{-2}}$ tends weakly to that of the multivariate normal
%
\begin{equation}
\NN_t :=\NN(0,t\frakq)
\end{equation}
for some positive semidefinite, nondegenerate $d\times d$ matrix $\frakq
=(q_{ij})$. As a consequence, if $\GG\subset C_0^\infty(\R^d)$ is a
family of functions that are uniformly equicontinuous and bounded, then
we have
%
\begin{equation}
\label{L2limit}
E_{\mu} \Bigl( \sup_{g\in\GG} | E_\kappa^0(g(\varepsilon X_{t\varepsilon
^{-2}}))-Eg(\NN_t) |^2 \Bigr) \mathop{\longrightarrow}_{\varepsilon\downarrow0} 0.
\end{equation}
Now, let $Q$ be the generator of the Brownian motion with mean zero and
covariance $\frakq$, that is, $Q$ is the operator in \eqref
{Q-operator}. We then have
%
\begin{equation}
\label{fQf}
(f,\texte^{ tQ}f)=\int\textd y\, f(y)E\bigl(f(y+\NN_t)\bigr).
\end{equation}
As $\LL_\kappa$ is the generator of the random walk $(X_t)$, we
similarly derive
%
\begin{eqnarray}\qquad
\label{fLf}
&&\varepsilon^{-2}(f_\varepsilon,\texte^{ t\varepsilon^{-2}\LL_\kappa}f_\varepsilon)
\nonumber\\[-8pt]\\[-8pt]
&&\qquad= \varepsilon^d \int_{[0,1]^d\times[0,1]^d} \textd z_1\,\textd z_2 \sum_{x\in
\Z^d}f(\varepsilon x+\varepsilon z_1)E_\kappa^x\bigl(f(\varepsilon z_2+\varepsilon
X_{t\varepsilon^{-2}})\bigr).\nonumber
\end{eqnarray}
We thus need to show that the right-hand side of \eqref{fLf} tends to
that of \eqref{fQf}. Note that if $f$ is supported in $[-M,M]^d$, then
the integral in \eqref{fQf} can be restricted to $y\in[-M,M]^d$ and the
sum over $x$ in \eqref{fLf} to, say, $|x|\le2M/\varepsilon$ (once
$\varepsilon\ll1$).

Substituting $y=\varepsilon x+\varepsilon z$ with $x\in\Z^d$, $|x|\le
2M/\varepsilon$ and $z\in[0,1]^d$ in \eqref{fQf} allows us to put both
terms on the same footing. Subtracting \eqref{fQf} from \eqref{fLf},
taking expectation with respect to $\mu$ and applying the
Cauchy--Schwarz inequality, we thus get
%
\begin{eqnarray}
\label{long}
&&E_\mu\Theta_\varepsilon(t)^2
\nonumber\\
&&\qquad\le\Vert f\Vert_2 \varepsilon^d \int_{[0,1]^d} \textd z \sum_{x\dvtx|x|\le
2M/\varepsilon}
E_{\mu} | E_\kappa^xf(\varepsilon z+\varepsilon X_{t\varepsilon
^{-2}})\\
&&\qquad\quad\hspace*{139.2pt}{}-Ef(\varepsilon x+\varepsilon z+\NN_t)
|^2.\nonumber
\end{eqnarray}
Using the translation invariance of $\mu$, we may replace $E_\kappa
^xf(\varepsilon z+\varepsilon X_{t\varepsilon^{-2}})$ by the expression $E_\kappa
^0f(\varepsilon x+\varepsilon z+\varepsilon X_{t\varepsilon^{-2}})$ inside the
expectation. For $f\in C_0^\infty(\R^d)$,
%
\begin{equation}
\GG:=\{f(\varepsilon x+\varepsilon z+\cdot)\dvtx|x|\le2M/\varepsilon, z\in
[0,1]^d\}
\end{equation}
is an equicontinuous family of uniformly bounded functions. Then, \eqref
{L2limit} tells us that the right-hand side of \eqref{long} tends to
zero as $\varepsilon\downarrow0$, which proves the desired claim \eqref{4.9}.
\end{pf}
\begin{pf*}{Proof of Proposition \ref{lemma-scaling}}
To extract \eqref{4.8} from \eqref{4.9}, we note that for any $f\in
C_0^\infty(\R^d)\cap\Hspace$,
%
\begin{equation}
\label{integral-repr}
(f,(-\LL_\kappa)^{-1}f)
=\int_0^\infty\textd t\, (f,\texte^{ t\LL_\kappa}f).
\end{equation}
Replacing $f$ by $f_\varepsilon$ and scaling $t$ by $\varepsilon^2$, we find
that
%
\begin{equation}
(f_\varepsilon,(-\LL_\kappa)^{-1}f_\varepsilon)=\int_0^\infty\textd t\, \varepsilon
^{-2}(f_\varepsilon,\texte^{t\varepsilon^{-2}\LL_\kappa} f_\varepsilon).
\end{equation}
By \eqref{4.9}, the function being integrated tends to $(f,\texte^{
tQ}f)$ in probability for each~$t$; the monotonicity in $t$ (and
continuity of the limit) ensures that the convergence is actually
uniform (in probability) on compact intervals. By Corollary \ref
{cor-tightness}, the integral can be truncated to a finite interval in
$L^1$-norm and similarly for the integral of the limit, which is finite
since $f$ is finite in the domain of $Q^{-1}$. It follows that
%
\begin{equation}
(f_\varepsilon,(-\LL_\kappa)^{-1}f_\varepsilon)
\mathop{\longrightarrow}_{\varepsilon\downarrow0} \int_0^\infty\textd t
\,(f,\texte^{ tQ}f)
=(f,-Q^{-1}f)
\end{equation}
in $\mu$-probability [and $L^1(\mu)$].
This is the desired conclusion \eqref{4.8}.
\end{pf*}
\begin{remark}
We note that to control the tail of the integral in \eqref
{integral-repr} in $d\ge3$, it suffices to invoke the diagonal heat
kernel estimate
%
\begin{equation}
E_\mu P_\kappa^0(X_t=x)\le\frac{c_1}{t^{d/2}},\qquad x\in\Z^d,
\end{equation}
which, in the elliptic case, is an immediate consequence of the mixing
theory for Markov chains based on isoperimetric inequalities. This is
sufficient because the finiteness of the Green function in $d\ge3$
permits us to define $(f,(-\LL_\kappa)^{-1}f)$, even for $f\ge0$. This
enables us to reduce the general case to positive $f$ by decomposing
the test function into a positive and a negative part and applying
%
\begin{equation}
E_\nu( f,\texte^{t\LL_\kappa}f)\le\Vert f\Vert_1^2\frac{c_1}{t^{d/2}},
\end{equation}
which is uniformly integrable when $d\ge3$.
However, to include $d=2$, we cannot disregard the cancellations due to
the vanishing of $\int f(x)\,\textd x$ and thus the stronger derivative
bound \eqref{HCUB} is necessary. A similar situation occurred in
Giacomin, Olla and Spohn \cite{GiacominOllaSpohn} where a stronger
Nash continuity estimate was required to include $d=2$.
\end{remark}

We are now ready to establish the main result of this paper.
\begin{pf*}{Proof of Theorem \ref{thm-main}}
Let $\mu$ be a translation-invariant, ergodic, gradient Gibbs measure
with zero tilt and let $\tilde\mu$ be its extension to $\R^{\B(\Z
^d)}\times\R^{\B(\Z^d)}$. We want to prove that $\phi(f_\varepsilon)$
tends weakly to a normal random variable with mean zero and variance
$(f,(-Q)^{-1}f)$. By Lemma \ref{lemma-cont}, it suffices to prove this
for $f\in C_0^\infty(\R^d)\cap\Hspace$.

By Lemma \ref{lemma-Gaussian}, we know that $\phi(f)$ is Gaussian
conditional on $\kappa$. The standard formula for any Gaussian random
variable $X$,
%
\begin{equation}
E(\texte^{\texti X})=\texte^{\texti E(X)-1/2\operatorname{Var}(X)},
\end{equation}
implies, via (\ref{mean-zero}) and (\ref{variance}), that
%
\begin{equation}
\label{4.9a}
E_{\tilde\mu}\bigl(\texte^{\texti\phi(f)}|\scrF\bigr)(\kappa)=\texte^{-
1/2(f,(-\LL_\kappa)^{-1}f)}.
\end{equation}
By Proposition \ref{lemma-scaling}, we have $(f_\varepsilon,(-\LL_\kappa
)^{-1}f_\varepsilon,)\to(f,(-Q)^{-1}f)$ in $\tilde\mu$-pro\-ba\-bi\-lity.
Since the right-hand side of \eqref{4.9a} is a bounded continuous
function of this inner product, \eqref{limit} follows by means of the
bounded convergence theorem.
\end{pf*}

\section{Potential theory for random conductance models}
\label{sec5}
The proof of the key Lemma \ref{lemma-harmonic} leads us to the study
of potential theory for operators depending on a random environment
that fall into the class of \textit{random conductance models}. A~good
deal of what is to follow exists explicitly, or implicitly, in the
literature. We have borrowed some of the notation from the paper of
Mathieu and Piatnitski \cite{MathieuPiatnitski}, although the
formalism draws on earlier works in homogenization theory; see, for
instance, the book by Jikov, Kozlov and Oleinik \cite{JKO}.
Notwithstanding, the content of Section \ref{sec5.2} appears to be new.

\subsection{Basic notions}
Consider a translation-invariant $\nu$ probability measure on $\Omega
:=\R_+^{\B(\Z^d)}$ (endowed with the product $\sigma$-algebra)
satisfying the ellipticity condition
%
\begin{equation}
\exists\varepsilon>0\dvtx
\nu\biggl(\varepsilon\le\kappa_b\le\frac1\varepsilon\biggr)=1,\qquad
b\in\B(\Z^d).
\end{equation}
Let $L^2(\nu)$ denote the closure of the set of all local functions in
the topology induced by the inner product
%
\begin{equation}
\label{bracket-product}
\langle h,g\rangle:=E_\nu(h(\kappa)g(\kappa)).
\end{equation}
Let $B:=\{\hate_1,\ldots,\hate_d\}$ denote the set of coordinate vectors
in $\Z^d$. The translations by the vectors in $B$ induce natural
unitary maps $T_1,\ldots,T_d$ on $L^2(\nu)$ defined via
%
\begin{equation}
(T_j h):=h\circ\tau_{\hate_j},\qquad j=1,\ldots,d.
\end{equation}
Apart from square integrable functions, we will also need to work with
vector fields, by which we will generally mean measurable functions
$u\dvtx\Omega\times B\to\R$ or $\Omega\times B\to\R^d$, depending on
the context. We will sometimes write $u_1,\ldots,u_d$ for $u(\cdot,\hate
_1),\ldots,u(\cdot,\hate_d)$---note that these may still be vector-valued.
\begin{remark}
While we index vector fields only by the \textit{positive} coordinate
vectors, in certain situations, it is convenient to have them also
defined for the negative coordinate directions via
%
\begin{equation}
\label{minus-b}
u(\kappa,-b):=-u(\tau_{-b}\kappa,b),\qquad b\in B.
\end{equation}
As we will see, this definition will automatically ensure that the
cycle condition (see Lemma \ref{lemma-cycle} below) holds for the
trivial cycles crossing only a single edge.
\end{remark}

Let $\Lvec^2(\nu)$ be the set of all vector fields with $(u,u)<\infty$,
where $(\cdot,\cdot)$ denotes the inner product
%
\begin{equation}
(u,v):=E_{\nu} \biggl( \sum_{b\in B}\kappa_b u(\kappa,b)\cdot v(\kappa,b) \biggr).
\end{equation}
Examples of such functions are the gradients $\nabla h$ of local
functions $h\in L^2(\nu)$ defined component-wise via the formula
%
\begin{equation}
(\nabla h)_j:=T_jh-h,\qquad j=1,\ldots,d.
\end{equation}
We denote by $L^2_{\nabla}(\nu)$ the closure of the set of gradients of
local functions in the topology induced by the above inner product.
\begin{lemma}
\label{lemma-cycle}
Let $u\in L^2_{\nabla}$. Then, $u$ satisfies the cycle condition
%
\begin{equation}
\label{cycle}
\sum_{j=0}^n u(\tau_{x_j}\kappa,x_{j+1}-x_j)=0
\end{equation}
for any finite (nearest-neighbor) cycle $(x_0,x_1,\ldots,x_n=x_0)$ on $\Z
^d$. In particular, there exists a shift-covariant function $\bar u\dvtx
\Omega\times\Z^d\to\R^d$ such that $u(\kappa,b)=\bar u(\kappa,b)$ for
every $b\in B$.
\end{lemma}
\begin{pf}
The cycle condition \eqref{cycle} holds trivially for all gradients of
local functions. Indeed, if $u=\nabla h$, then, in light of \eqref
{minus-b}, we have
%
\begin{equation}\qquad
u(\tau_{x_j}\kappa,x_{j+1}-x_j)=(\nabla h)_{x_{j+1}-x_j}(\tau
_{x_j}\kappa)=h\circ\tau_{x_{j+1}}(\kappa)-h\circ\tau_{x_j}(\kappa).
\end{equation}
A corresponding limit extends this to all of $L^2_{\nabla}$. To define
$\bar u(\cdot,x)$, we integrate properly shifted values of $u$ along a
path from zero to $x$; the cycle condition guarantees that the result
is independent of the choice of path and that $\bar u$ is shift-covariant.
\end{pf}

We will henceforth use the convention of writing $\bar u$ for the
extension of a shift-covariant vector field $u\in\Lvec^2$ to a function
on $\Z^d$.
Notice that the shift $T_j$ extends naturally via
%
\begin{equation}
T_j\bar u(\kappa,x):=\bar u(\tau_{\hate_j}\kappa,x)=\overline{(T_j
u)}(\kappa,x).
\end{equation}
Next, let us characterize the functions in $(L^2_{\nabla})^\perp$.
\begin{lemma}
\label{lemma-divergence}
For $u\in\Lvec^2(\nu)$, let $\LL u$ be the function in $L^2(\nu)$
defined by
%
\begin{equation}
(\LL u)(\kappa):=\sum_{b\in B}[ \kappa_b u(\kappa,b)-(\tau_{-b}\kappa
)_{b} u(\tau_{-b}\kappa,b)],
\end{equation}
where $-b$ is the coordinate vector opposite to $b$. We then have
%
\begin{equation}
u\in(L^2_{\nabla})^\perp\quad\Leftrightarrow\quad\LL u=0,\qquad \nu\mbox{-a.s.}
\end{equation}
If $u$ satisfies the cycle condition and $\bar u$ is its extension,
then $\LL u(\tau_x\kappa)=\LL_\kappa\bar u(\kappa,x)$.
\end{lemma}
\begin{pf}
These are direct consequences of the definitions, the translation
invariance of $\nu$ and a simple calculation.
\end{pf}

Note that $\LL u$ plays the role of the \textit{divergence}---that is,
the total flow out of a given vertex---of vector field $u$. However, we
prefer to denote it by $\LL$ to emphasize its connection with the
operator $\LL_\kappa$.

\subsection{Uniqueness of harmonic embedding}
\label{sec5.2}
Clearly, all $u\in L^2_{\nabla}$ are shift-covariant and have zero
mean. A question which naturally arises is whether every
shift-covariant zero-mean $u$ is in $L^2_{\nabla}$. (Note that this is
analogous to asking whe\-ther every closed differential form is exact.)
Our answer to this is in the affirmative.
\begin{theorem}
\label{thm-complement}
Suppose $\nu$ is ergodic. Then, every $u\in\Lvec^2(\nu)$ which obeys
the cycle condition \eqref{cycle} and $E_\nu u=0$ satisfies $u\in
L^2_{\nabla}$.
\end{theorem}

Again, recall that \eqref{cycle} and zero expectation are necessary for
$u\in L^2_{\nabla}$. The above implies that these conditions are also
sufficient. To prove the theorem, we will need the following lemma.
\begin{lemma}
\label{lemma-projection}
Let $P_j$ denote the orthogonal projection onto $\operatorname
{Ker}(1-T_j)$ in $\Lvec^2(\nu)$. If $\nu$ is ergodic and $u\in\Lvec
^2(\nu)$ satisfies \eqref{cycle}, then $P_j u=E_\nu(P_j u)$, $\nu
$-almost surely.
\end{lemma}
\begin{pf}
Fix $u\in\Lvec^2(\nu)$ that obeys \eqref{cycle} and let $\bar u$ be
the corresponding shift-covariant function. We will prove the claim
only for the component $u_1=u(\cdot,\hate_1)$; the other cases follow
analogously. By translation covariance and the $L^2$ ergodic theorem,
we have
%
\begin{equation}
\label{on-axis}
\frac{\bar u(\cdot,n\hate_1)}n=\frac1n\sum_{k=0}^{n-1}T_1^ku_1 \mathop
{\longrightarrow}_{n\to\infty} P_1u_1\qquad
\mbox{in }L^2(\nu).
\end{equation}
If $\nu$ were separately ergodic (i.e., ergodic with respect to $T_1$
alone), then the claim would immediately follow by the fact that every
$T_1$-invariant function must be constant. To make up for the potential
lack of separate ergodicity, we note that translation covariance of $u$
and the fact that $\bar u$ obeys the cycle conditions together yield
%
\begin{eqnarray}
T_j\bar u(\cdot,n\hate_1)&
=&\bar u(\cdot, n\hate_1+\hate_j)-\bar u(\cdot,\hate_j)
\nonumber\\[-8pt]\\[-8pt]
&=&\bar
u(\cdot,n\hate_1)-u(\cdot,\hate_j)+T_{n\hate_1}u(\cdot,\hate_j).\nonumber
\end{eqnarray}
It follows that $T_j \frac1n \bar u(\cdot,n\hate_1)$ also converges to
$P_1u_1$ (in $L^2$) and so, by the continuity of $T_j$,
%
\begin{equation}
T_j P_1u_1=P_1u_1.
\end{equation}
Hence, $P_1u$ is invariant under all shifts and is therefore constant
$\nu$-a.s.
\end{pf}
\begin{pf*}{Proof of Theorem \ref{thm-complement}}
Suppose $\nu$ is ergodic and let $u\in\Lvec^2(\nu)$ obey \eqref{cycle}
and $E_\nu u=0$. The boundedness of the $\kappa_b$'s away from zero and
infinity ensures that $u\in\Lvec^2(\nu)$ if and only if all of its
components are in $L^2(\nu)$. Our task is to construct functions
$h_\varepsilon\in L^2(\nu)$ such that $\nabla h_\varepsilon\to u$ in $\Lvec
^2(\nu)$. We define
%
\begin{equation}
\label{hepsilon}
h_\varepsilon:=-\sum_{n\ge0}\frac{T_1^nu_1}{(1+\varepsilon)^{n+1}}
\end{equation}
and note that this is the unique solution of the equation $(1+\varepsilon
-T_1)h_\varepsilon=-u_1$. This observation implies
that
%
\begin{equation}
(1-T_1)h_\varepsilon=-u_1-\varepsilon h_\varepsilon
\end{equation}
and so the first component of $\nabla h_\varepsilon$ converges to that of
$u$, provided $\varepsilon h_\varepsilon\to0$ in $L^2(\nu)$. To see what
happens with the other components of $\nabla h_\varepsilon$, we note that
the cycle condition \eqref{cycle} translates into
%
\begin{equation}
(1-T_j)u_1=(1-T_1)u_j.
\end{equation}
Applying this to the definition of $h_\varepsilon$, we conclude
that
%
\begin{equation}
(1-T_j)h_\varepsilon=-u_j-\varepsilon\tilde h_\varepsilon,
\end{equation}
where $\tilde h_\varepsilon$ is defined as $h_\varepsilon$, but with $u_1$
replaced by $u_j$. Again, it suffices to show that $\varepsilon\tilde
h_\varepsilon\to0$ in $L^2(\nu)$, which will boil down to the same
argument as for $j=1$.

To prove that $\varepsilon h_\varepsilon\to0$, we note that, for the inner
product in \eqref{bracket-product},
%
\begin{equation}
\langle h_\varepsilon,h_\varepsilon\rangle=\sum_{n\ge0}\frac{n+1}{(1+\varepsilon
)^{n+2}}\langle u_1,T_1^nu_1\rangle.
\end{equation}
Introducing the notation $A_nu$ for the average
%
\begin{equation}
A_nu:=\frac1n\sum_{k=0}^{n-1}T_1^ku,
\end{equation}
inserting this into the above sum and reordering the terms, we get
%
\begin{equation}
\label{hh}
\langle h_\varepsilon,h_\varepsilon\rangle=
\frac{\langle u_1,u_1\rangle}{(1+\varepsilon)^2}+
\sum_{n\ge1}\frac{n(\varepsilon n-1)}{(1+\varepsilon)^{n+2}} \langle
u_1,A_nu_1\rangle.
\end{equation}
By Lemma \ref{lemma-projection}, the $L^2$ ergodic theorem and the fact
that $u$ has zero expectation in~$\nu$,
we have
%
\begin{equation}
\label{Au1-constant}
A_nu_1 \mathop{\longrightarrow}_{n\to\infty} P_1u_1=E_\nu u_1=0 \qquad\mbox{in
}L^2(\nu)
\end{equation}
and so $\langle u_1,A_nu_1\rangle\to0$ as $n\to\infty$. A
straightforward estimate now shows that the sum in \eqref{hh} is
$o(\varepsilon^{-2})$ and so $\varepsilon^2\langle h_\varepsilon,h_\varepsilon
\rangle\to0$, as desired.
\end{pf*}

Not every function in $\Lvec^2(\nu)$ necessarily belongs to $L^2_{\nabla
}(\nu)$. A prime example is the \textit{position} vector field $x(\kappa
,b)=b$. Indeed, let $\chi\dvtx\R^{\B(\Z^d)}\times B\to\R^d$ be the projection
%
\begin{equation}
\chi:=-\operatorname{proj}_{L^2_{\nabla}(\nu)}x.
\end{equation}
Since $\chi\in L^2_{\nabla}(\nu)$, it satisfies \eqref{cycle} and we
may extend it to a function $\bar\chi$ mapping $\R^{\B(\Z^d)}\times\Z
^d\to\R^d$ by setting $\bar\chi(\cdot,0)=0$ and integrating the
gradients along oriented paths. Lemma \ref{lemma-divergence} implies
that $\bar x+\bar\chi$ is harmonic, $\LL_\kappa(\bar x+\bar\chi)=0$.
Moreover, Lemma \ref{lemma-projection} and \eqref{on-axis} show that
$\bar\chi(\cdot,n\hate_j)/n\to0$ and so $x+\chi\ne0$. It follows that
$x\notin L^2_{\nabla}$ and so $L^2_{\nabla}\ne\Lvec^2$.

The function $\bar\chi$ is generally referred to as the \textit
{corrector} because it corrects for the nonharmoni\-city of the
position function. The corrector can be defined by appealing to
spectral theory (Kipnis and Varadhan \cite{KipnisVaradhan}, also
Berger and Biskup \cite{BergerBiskup}); the above ``projection''
definition is inspired by those in Giacomin, Olla and Spohn
\cite{GiacominOllaSpohn} and Mathieu and Piatnitski \cite{MathieuPiatnitski}.
As a side remark, we note that the function $\bar x+\bar\chi$ actually
allows us to characterize the space of all square integrable
shift-covariant functions.
\begin{corollary}
\label{cor-characterize}
Suppose $\nu$ is ergodic. Then, every shift-covariant $\R$-valued $u\in
(L^2_{\nabla})^\perp$ can be obtained from $x+\chi$, where $x$ is the
position function and $-\chi$ is its orthogonal projection onto
$L^2_{\nabla}$, by means of the linear transformation
%
\begin{equation}
u_j(\kappa)=\sum_{k=1}^d [\hate_k\cdot(x_j+\chi_j)]E_\nu u_k,\qquad j=1,\ldots,d.
\end{equation}
Here, $u_1,\ldots,u_d$ stand for $u(\cdot,\hate_1),\ldots,u(\cdot,\hate
_d)$ and similarly for ($\R^d$-valued objects) $x_1,\ldots,x_k$ and $\chi
_1,\ldots,\chi_d$. In particular, for $\R$-valued vector fields, we have
%
\begin{equation}
\label{L2-decompose}
\{u\in\Lvec^2(\nu)\dvtx\mbox{shift-covariant}\}=L^2_{\nabla} \oplus\{
\lambda\cdot(x+\chi)\dvtx\lambda\in\R^d\}.
\end{equation}
\end{corollary}
\begin{remark}
The observation \eqref{L2-decompose}---which, in the language borrowed
from differential geometry, implies that linear transforms of $x+\chi$
are the only closed forms that are not exact---was previously
discovered in the context of multicolored exclusion processes (Quastel
\cite{Quastel}, Theorem 9.1). However, while our proof is based only on
a soft, Poisson equation-based argument \twoeqref
{hepsilon}{Au1-constant}, that of \cite{Quastel} requires an explicit
bound on the spectral gap for the corresponding dynamics. The reason
may be the reliance of \cite{Quastel} on the \textit{spatial} ergodic
theorem, which, naturally, leads to bounds involving the Poincar\'e
inequality and/or spectral gap.
\end{remark}
\begin{pf*}{Proof of Corollary \ref{cor-characterize}}
This is a simple consequence of Theorem \ref
{thm-complement}. Let $w=(w_1,\ldots,w_d)$ be the $d$-component vector
field whose ($\R$-valued) components are defined by
%
\begin{equation}
w_j:=\sum_{k=1}^d [\hate_k\cdot(x_j+\chi_j)]E_\nu u_k,\qquad j=1,\ldots,d.
\end{equation}
We will show that $w=u$. First, both $u$ and $w$ obey the cycle
condition, so $u-w$ also does. As $\chi\in L^2_{\nabla}(\nu)$ implies
$E_\nu\chi=0$, the fact that $E_\nu x_k=x_k=\hate_k$ shows that $E_\nu
w_j=E_\nu u_j$, that is, $E_\nu(u-w)=0$. Theorem \ref
{thm-complement} implies that $u-w\in L^2_{\nabla}$. On the other
hand, $\LL(x+\chi)=0$ implies that $\LL w=0$ and so, by Lemma \ref
{lemma-divergence}, $w\in(L^2_{\nabla})^\perp$. Thus, $u-w\in
(L^2_{\nabla})^\perp$ also holds. It follows that $u=w$, as claimed.

As for \eqref{L2-decompose}, the argument we just used ensures that a
shift-covariant field $u\in\Lvec^2(\nu)$ can be written as $\lambda\cdot
(x+\chi)$, where $\lambda$ is a vector with components $\lambda
_j:=E_\nu u_j$, plus a shift-covariant vector field with zero
expectation. By Theorem~\ref{thm-complement}, the latter is in
$L^2_\nabla$.
\end{pf*}

We still have to supply the proof of Lemma \ref{lemma-harmonic}.
\begin{pf*}{Proof of Lemma \ref{lemma-harmonic}}
Since $g$ has square integrable components, we have $g\in\Lvec^2(\nu)$.
As $g$ is shift-covariant and has zero expectation, Theorem \ref
{thm-complement} implies that $g\in L^2_{\nabla}$. However, $g$ is also
harmonic and so, in turn, Lemma \ref{lemma-divergence} forces $g\in
(L^2_{\nabla})^\perp$. Thus, $g=0$, as desired.
\end{pf*}

Recall that $\bar\chi$ is the extension of $\chi$ subject to the
condition $\bar\chi(\cdot,0)=0$. With this object at hand, we may even
remove the restriction to zero slope in Lem\-ma~\ref{lemma-Gaussian}.
\begin{corollary}
\label{cor-tilt}
Let $\mu$ be a translation-invariant, ergodic, gradient Gibbs measure
for the potential \eqref{V-repr} and let $\tilde\mu$ be its extension
to $\R^{\B(\Z^d)}\times\R^{\B(\Z^d)}$. Let $t\in\R^d$ be the tilt of
$\mu$ and let $\scrF:=\sigma(\{\kappa_b\dvtx b\in\B(\Z^d)\})$. Then,
$\tilde\mu(-|\scrF)(\kappa)$ is Gaussian with mean
%
\begin{equation}
\label{mean}
E_{\tilde\mu}(\phi_x-\phi_0|\scrF)(\kappa)=t\cdot[x+\bar\chi(\kappa,x)]
\end{equation}
and covariance given by $(-\LL_\kappa)^{-1}$.
\end{corollary}
\begin{pf}
Let $\bar u$ denote the extension of the shift-covariant vector field
$u(\kappa, b):=E_{\tilde\mu}(\eta_{0,b}|\scrF)(\kappa)$, $b\in B$.
Clearly, we have
%
\begin{equation}
\bar u(\kappa,x)=E_{\tilde\mu}(\phi_x-\phi_0|\scrF)(\kappa),\qquad x\in\Z^d.
\end{equation}
Inspecting the proof of Lemma \ref{lemma-Gaussian}, all formulae carry
over, except \eqref{u-exp}, which becomes
%
\begin{equation}
\label{u-exp2}
E_{\tilde\mu}(\bar u(\cdot,x))=t\cdot x,\qquad x\in\Z^d.
\end{equation}
Thus, $\bar u-t\cdot(\bar x+\bar\chi)$ is harmonic, shift-covariant and
of zero mean, so $\bar u=t\cdot(\bar x+\bar\chi)$ by Lemma \ref
{lemma-divergence} and Theorem \ref{thm-complement}.
\end{pf}

\section{Discussion and open problems}
\label{sec7}
The proofs in the present note often rely on the fact that the random
walk with generator $\LL_\kappa$ is uniformly elliptic. This enters via
the assumption that $\varrho$, defining the potential $V$, has support
bounded away from zero. While we believe that the general picture
carries over, even if we let the support extend all the way to zero, a
number of steps in the proof become quite subtle. For instance, the
pointwise heat kernel asymptotic for this walk may take a radically
different form (Berger et al. \cite{BBHK}) and it
is not known under what conditions on the environment the walk scales
to Brownian motion.
Progress in this direction for i.i.d. (or i.i.d. dominated)
environments has been made only recently (Mathieu \cite{MathieuCLT},
Biskup and Prescott \cite{BiskupPrescott}, Barlow and
Deuschel~\cite{BarlowDeuschel}).

Another interesting open question concerns the annoying restriction to
\textit{zero-tilt} gradient Gibbs measures. As is seen from Corollary \ref
{cor-tilt}, once the tilt is nonzero, \eqref{4.9a} has to be modified to
%
\begin{equation}
E_{\tilde\mu}\bigl(\texte^{\texti[\phi(f)-E_{\tilde\mu}(\phi(f))]}|\scrF
\bigr)(\kappa)=\texte^{-\texti(t\cdot\bar\chi)(f)-1/2(f,(-\LL_\kappa)^{-1}f)}.
\end{equation}
When we plug in $f_\varepsilon$ for $f$, the second term in the exponent
still converges to $(f,(-Q)^{-1}f)$ in probability as $\varepsilon
\downarrow0$. To show convergence of $\phi_\varepsilon-E_{\tilde\mu}\phi
_\varepsilon$ to Gaussian free field, we thus have to show that $(t\cdot
\bar\chi)(f_\varepsilon)$ converges in law to a normal random variable. If
$V$ is strictly convex and of the form \eqref{V-repr}, then we know
this to be true by the main result of Giacomin, Olla and Spohn
\cite{GiacominOllaSpohn}. However, the general case seems to be a hard
open problem (cf. \cite{BergerBiskup}, Conjecture 5). Already in
$d=1$, where the corrector can be written explicitly as the sum
%
\begin{equation}
\bar\chi(\kappa,x)=\frac1C\sum_{n=0}^{x-1} \biggl( \frac1{\kappa_{n,n+1}}-C \biggr),
\end{equation}
with $C:=E_{\tilde\mu}(1/\kappa_b)$, the desired result requires that
the resistivities $1/\kappa_b$ exhibit a specific correlation decay
\cite{Durrett}, Theorem 7.6. Very little information concerning this is
known in $d\ge2$; see Mourrat \cite{Mourrat} for some recent progress
in this direction.

Our final remark concerns a generalization to potentials $V$, where
\eqref{V-repr} has been modified to
%
\begin{equation}
V(\eta):= -\log\int\varrho(\textd\kappa) \texte^{-W_\kappa(\eta)},
\end{equation}
with $(W_\kappa)$ denoting a family of strictly convex, even,
measurable functions with uniformly superlinear growth at $\pm\infty$
and a uniform lower bound. In this case, we may still consider the
extended gradient Gibbs measures; however, the conditional law given
the $\kappa$'s is no longer Gaussian. Notwithstanding, given $\kappa$,
the Helffer--Sj\"ostrand
representation still applies and leads to a random walk in a dynamic
random environment. Its annealed central limit theorem would imply the
Gaussian scaling limit for the $\phi$-field.

\section*{Acknowledgments}

This project was initiated at the Kavli Institute for
Theoretical Physics during the program ``Stochastic Geometry and Field
Theory: From Growth Phenomena to Disordered Systems.'' The KITP program
was supported by the NSF Grant PHY99-07949.
We wish to thank Jean-Dominique Deuschel for a number of interesting
suggestions, and Scott Sheffield and Jeremy Quastel for discussions on
related subjects. We are also grateful to the anonymous referees for
careful reading of the first version of this paper and for a number of
constructive suggestions.


%
\printaddresses

\end{document}